\documentclass{volcon}
\usepackage{amssymb,epsfig}

\begin{document}

\begin{frontmatter}

\title{A Non-Existence Result for\\ Hamiltonian Integrators}

\author{P. F. Tupper}
\ead{tupper@math.mcgill.ca}
\ead[url]{www.math.mcgill.ca/\~{ }tupper}

\address{Department of Mathematics and Statistics, McGill University,\\
Montr\'{e}al QC, H3A 2K6 Canada.}


\begin{abstract}
We consider the numerical simulation of Hamiltonian systems of ordinary differential equations.  Two features of Hamiltonian systems are that energy is conserved along trajectories and phase space volume is preserved by the flow.
We want to determine if there are integration schemes that preserve these two properties for all Hamiltonian systems, or at least for all systems in a wide class.
This paper provides provides a negative result in the case of two dimensional (one degree of freedom) Hamiltonian systems, for which phase space volume is identical to area.
Our main theorem shows that there are no computationally reasonable numerical integrators for which all Hamiltonian systems of one degree of freedom can be integrated while conserving both area and energy.  Before proving this result we define what we mean by a computationally reasonable integrator.  We then consider what obstructions this result places on the existence of volume- and energy-conserving integrators for Hamiltonian systems with an arbitrary number of degrees of freedom.
\end{abstract}

\end{frontmatter}

\section{Introduction}

We consider a system of Hamiltonian differential equations on $\mathbb{R}^{2n}$ 
\begin{equation} \label{eq:Ham0}
\frac{ dq}{dt} = \frac{\partial H}{\partial p}, \ \ \ \ \frac{dp}{dt} = -\frac{\partial H}{\partial q},
\end{equation}
defined by the Hamiltonian function $H: \mathbb{R}^{2n} \rightarrow
\mathbb{R}$.
We denote the time $t$ flow map of these equations by $S^t$. 
The flow of these differential equations has two important features.  The first is that $H$ is conserved along trajectories.  That is,
\(
H(S^t (q,p)) = H(q,p),
\)
for all $t\in \mathbb{R}$ and all $(q,p) \in \mathbb{R}^{2n}$. 
The second is that phase space volume is conserved by the flow:
if $A \subset \mathbb{R}^{2n}$ is a bounded open set, then 
\(
\mathrm{vol} (S^t A) = \mathrm{vol}(A)
\)
for all $t$.
The latter property is a consequence of the symplecticity of the flow.  (See, for example, \cite{hairer}).

For certain molecular dynamics applications, the ideal numerical integrator
would retain these two properties of the flow \cite{tupper}.
That is, if the integrator with time step $\Delta t$ defines a map $\Phi_{\Delta t}$ on
$\mathbb{R}^{2n}$, we would like both $H(\Phi_{\Delta t}(q,p)) = H(q,p)$ for all $(q,p)$
(energy conservation)
and $\mbox{vol}(\Phi_{\Delta t}( A) ) = \mbox{vol}(A)$ for all bounded open
subsets $A$ of $\mathbb{R}^{2n}$(volume conservation). 

Symplectic integrators such as the implicit midpoint rule conserve
volume exactly for all Hamiltonian systems, with any number of degrees of freedom, 
but do not conserve energy \cite{hairer}.   It has 
already been shown that in certain circumstances it is unreasonable to
expect that symplectic integrators which also conserve energy exist
\cite{zhong}.  However,  volume-conservation is a weaker property than
symplecticity for Hamiltonian systems of more than one degree of
freedom ($n>1$). 
Thus it seems plausible that there is a consistent integration scheme which for any Hamiltonian function $H$ and time step $\Delta t$ yields a map $\Phi_{\Delta t}$ which conserves both volume and energy.

In this article we will argue that there is no such integration scheme.
We do this by showing that no numerical integrator
 is able to  integrate \emph{all} Hamiltonian systems in $\mathbb{R}^2$ while simultaneously  conserving energy and
phase space volume.   For special Hamiltonian systems in $\mathbb{R}^2$ such an integrator is possible.  (See \cite{quispel} for examples of such systems of any number of degrees of freedom.)
However, our theorem states that for any energy-conserving numerical
integrator from a very broad class, there will be at least one Hamiltonian
system on $\mathbb{R}^2$ (and in fact very many) for which it does not conserve
phase space volume. 
 
Before stating and proving our main theorem, in 
Section~\ref{sec:integrator} we will define what we mean by an integrator.   From the class of integrators we then define the
  \emph{computationally reasonable} integrators.  This will include
  any 
  explicit or implicit formula for defining a new state $(q^{k+1},p^{k+1})$ as a
  function of a state $(q^k,p^k)$ and a timestep $\Delta t$ such that the only
  information used about $H$ is its value and the value of its derivatives at a finite
  number of points $(r^j,s^j) \in \mathbb{R}^{2n}$, $j=1,\ldots,m$.  In Section~\ref{sec:main} we will state and prove our main result in terms of this definition.  
  We imagine that a numerical
  analyst has devised a computationally reasonable numerical integrator that is
  energy-conserving for all Hamiltonian systems in $\mathbb{R}^2$.  We apply this
  integrator to the system with Hamiltonian
  function $H(q,p)= p^2/2$ in $\mathbb{R}^2$.  
  For two different inputs  we observe at which points $(r^j,s^j), j=1,\ldots,m$ the integrator depends on the function $H$ and its derivatives.  
Using this information, we construct another Hamiltonian function $\widetilde{H}(q,p)=p^2/2+V(q)$ which is arbitrarily close to $H$ and has the following property: the numerical integrator cannot be volume-conserving for both $H$ and $\widetilde{H}$.  
  
 The main result demonstrates that there is no
 integration scheme that conserves volume and energy for arbitrary Hamiltonian systems of any number of degrees of freedom.
 However, we cannot conclude that there is no scheme
 that is energy- and volume-conserving for Hamiltonian systems in a particular
 number of dimensions $2n$,
 $n>1$.  This question is open. To partially address this issue, in Section~\ref{sec:multi}
 we will show that under some
 reasonable--- but not essential--- conditions on an integrator, the
 problem reduces to that of the $n=1$ case.  Thus, there is no energy-
 and volume-conserving integrator for fixed $n$ that satisfies these
 additional assumptions. 

We now explain the relation between the results here and those of the well-known paper of Zhong and Marsden \cite{zhong}.  In that paper the authors consider a Hamiltonian system for which there are no invariants except energy.  They show that any integrator that is both symplectic and energy-conserving actually computes exact trajectories of the system up to a time reparametrization.  
From this result they conclude that energy-conserving symplectic integration is not possible in general, since presumably the set of Hamiltonians for which one could compute trajectories exactly up to a time reparametrization are very small.
Though this argument is very plausible, it leaves open two questions:
\begin{enumerate}
\item  How do we make precise the idea of something not being possible for any numerical integrator?
\item Is it possible to perform energy-conserving and symplectic integration for general Hamiltonians when we restrict ourselves to the $\mathbb{R}^2$ case?
\end{enumerate}
The importance of the latter question is that if it were possible, then volume- and energy-conserving integration would be possible for general Hamiltonian systems via $J$-splitting \cite{quispel,faou}.

We provide an answer to the first question in Section~\ref{sec:integrator} with the definition of a computationally reasonable integrator. 
As for the second question,
since volume-conservation and symplecticity are identical in $\mathbb{R}^2$, Zhong and Marsden's result shows that in $\mathbb{R}^2$ volume- and energy-conserving integration is equivalent to solving the original system exactly up to a time-reparametrization.  (This result is stated and proved for this case as Lemma~\ref{lem:reparam} in Section~\ref{sec:main}.)  
The main result of our paper answers the second question in the negative by  showing that it is not possible for general Hamiltonian systems in $\mathbb{R}^2$ using computationally reasonable integrators.

  Before we begin we discuss some interesting related work.  Even though energy and volume conserving integrators may not exist, the paper \cite{faou} does the next best thing. 
 There the authors show
  how to approximate any Hamiltonian function arbitrarily well by a
  special piece-wise smooth function of a form described by \cite{quispel} whose trajectories can be
  integrated while conserving volume and energy.  The original
  Hamiltonian function is not conserved.   However, unlike for
  standard symplectic methods, a modified Hamiltonian function close to the original is conserved exactly \cite{faou} for all time.

\section{What is a numerical integrator?}
\label{sec:integrator}

An integrator for a Hamiltonian system of ordinary differential equations 
takes a Hamiltonian $H$, a step length $\Delta t$, and an initial value
$(q^k,p^k)\in \mathbb{R}^{2n}$, and produces a value $(q^{k+1},p^{k+1}) \in \mathbb{R}^{2n}$.    Typically,
$(q^{k+1},p^{k+1})$ depends 
on $H$ through components of $\nabla H$ and perhaps $H$ itself. If
$(q^k,p^k)$  is an approximation to $(q(k \Delta t),p(k \Delta t))$, 
we take  $(q^{k+1},p^{k+1})$ to
be an approximation to $(q((k+1)\Delta t),p((k+1)\Delta t))$.
\begin{defn} \label{defn:integrator}
An integrator $\Phi$ is a
function that takes arguments $H : \mathbb{R}^{2n} \rightarrow
\mathbb{R}$, $\Delta t \in [0,\infty)$, $(q^k,p^k) \in \mathbb{R}^{2n}$ and
 either  returns $(q^{k+1},p^{k+1}) \in \mathbb{R}^{2n}$ or is not defined.
We write
\[
(q^{k+1},p^{k+1}) = \Phi( H, (q^k,p^k),\Delta t).
\]
\end{defn}
We have allowed the integrator to not be defined for certain input
values.  This is often the case for implicit integrators when the
vector field is insufficiently smooth or the time step is too large.

In
order to get meaningful constraints on what is computationally
reasonable, we cannot let arbitrary maps be included in the class
of algorithms we study.  After all, the exact flow map has all the qualitative features one could want, but it is not feasible to compute it (even to
machine precision) for most applications.
We would like our definition to be broad enough to include most existing numerical integrators.

Informally, we say an integrator $\Phi$ is
  \emph{computationally reasonable} if 
for each $(q^k,p^k)$ and $\Delta t$, $(q^{k+1},p^{k+1})=\Phi(H,(q^k,p^k),
\Delta t)$ depends 
on $H$ only through its value and the value of its derivatives at a
finite number of points $(r^j,s^j) \in \mathbb{R}^{2n}$, $j=1,\ldots,n$. 
In the following formal definition we use multi-index notation to define higher-order derivatives: for $H : \mathbb{R}^{2n} \rightarrow \mathbb{R}$ and $\alpha \in \mathbb{N}^{2n}_0$
we let 
\[
\partial_\alpha H := 
\frac{\partial^{\alpha_1}}{\partial q_1^{\alpha_1}}
\cdots
\frac{\partial^{\alpha_n}}{\partial q_n^{\alpha_n}}
\frac{\partial^{\alpha_{n+1}}}{\partial p_1^{\alpha_{n+1}}}
\cdots
\frac{\partial^{\alpha_{2n}}}{\partial p_n^{\alpha_{2n}}} H
\]

\begin{defn} \label{defn:reasonable}  An integrator $\Phi$ is
  \emph{computationally reasonable} if 
for each $H: \mathbb{R}^{2n} \rightarrow \mathbb{R}$, $(q^k,p^k) \in \mathbb{R}^{2n}$ and $\Delta t \geq
0$ there exists 
\begin{enumerate}
\item $m \in \mathbb{N}_0$,
\item $(r^j,s^j) \in \mathbb{R}^{2n}, j=1, \ldots, m$,
\item $\alpha^j \in \mathbb{N}_0^{2n}, j=1, \ldots, m$,
\end{enumerate}
such that for any function $\widetilde{H}: \mathbb{R}^{2n} \rightarrow \mathbb{R}$ that satisfies
\[
\partial_{\alpha^j} \widetilde{H} (r^j) = \partial_{\alpha^j} H(r^j), \ \ \ j=1,\ldots,m
\]
either $\Phi(H,(q^k,p^k),\Delta t) =
\Phi(\widetilde{H},(q^k,p^k),\Delta t)$ or both are not defined. 
\end{defn}
In the remainder of this section we discuss examples of the
class of computationally reasonable integrators.   

First note that all explicit methods fit into this class.  We informally define
an integrator to be explicit if it can be implemented by an algorithm
that terminates in a finite number of steps using function evaluations
of $H$ or its derivatives, arithmetic operations, and logical
operations.  This includes all of the explicit Runge-Kutta and
partitioned Runge-Kutta methods, for example.  It also includes
integrators that collect information adaptively to perform a step, 
such as the Bulirsch-Stoer method \cite{bs}.
Taylor series methods are also included. 

The class of computationally reasonable integrators also includes implicit methods, such as the implicit Runge-Kutta methods.  Here we define an implicit method to be one where $(q^{k+1},p^{k+1})$ is specified by requiring it to be the solution to a nonlinear system of equations in $H$ and its derivatives.
Implicit algorithms
cannot in general be implemented exactly in 
a finite number of steps, but they still fit into the framework of Definition~\ref{defn:reasonable}.  To see this, note that even though solving a system of nonlinear equations exactly typically requires looking at $H$ and its derivatives at an infinite number of points (while performing the Newton Iteration, for example), determining if we have a solution to a nonlinear system of equations only requires examining a finite number of points.  So if $(q^{k+1},p^{k+1})$ solves the equations for $H$, it will still solve the equation for any $\widetilde{H}$ which is identical to $H$ at the points.
Similarly, step-and-project methods are included in this class.  
These are methods that consist of one step of a simpler method
followed by a projection onto a manifold \cite[IV.4]{hairer}.  

What integrators do not satisfy Definition~\ref{defn:reasonable}?
Integrators that require the exact computation of integrals of $H$ or its derivatives do not.   Computing the integral of a general function requires knowing its value at an infinite set of points on the domain of integration.  Unlike in the case of solving nonlinear equations, if we are given a value for the integral, there is no finite set of points at which we can examine the function to verify that the value is correct.
 Of course, it is possible to define an numerical integrator that uses integrals of functions, but the actual computation of these integrals for general Hamiltonian systems would require numerical quadrature.  This in turn would require sampling the function at a finite number of points, and introducing truncation error.  The new method with this additional truncation error does form a computationally reasonable integrator, while the original method with the exact integral does not.

Finally, we note that multistep methods are not even integrators according to Definition~\ref{defn:integrator}.  We believe our framework could be extended to multistep method but we do not do so here.

\section{Main Theorem}
\label{sec:main}

To prove our main result Theorem~\ref{thm:main} we use the following lemma.  It shows that for a Hamiltonian system in $\mathbb{R}^{2}$, the map defined by an energy and area preserving integrator is just a time-reparametrization of the flow map.   As discussed in the introduction, this is essentially Zhong and Marsden's result \cite{zhong} in the two dimensional case.  The only addition is that we show that the time-reparametrization is just a constant rescaling of time where locally the constant does not depend on energy.

\begin{lem} \label{lem:reparam}
Let $H:\mathbb{R}^2 \rightarrow \mathbb{R}$ be a smooth function and
$S^t$ its induced Hamiltonian flow map.   Let $E \in \mathbb{R}$ be a particular energy.  Let $U$ be an open set whose intersection with  $\{ (q,p) | H(q,p) = E \}$ is a simple curve $\Gamma$.  Suppose that $\nabla H \neq 0$ on $\Gamma$.
Let $\Phi_{\Delta t}$, $\Delta t>0$ be a
continuous area-conserving map defined on $U$ that
conserves $H$. 
Then there is a constant  $c$ such that $\Phi_{\Delta t}(q,p) = S^{c \Delta t}(q,p)$ for all 
$(q,p) \in \Gamma$.
\end{lem}

{\bf Proof: }
On $U$
we can define canonical action-angle coordinates
$(\theta, \phi)$ in which the  Hamiltonian function is  $H(\theta,\phi) =
H(\theta)$.  The flow map is then 
\(
S^t(\theta,\phi) = (\theta,\phi + t H'(\theta) ),
\)
where $H'(\theta) \neq 0$, since we still have  $[H'(\theta), 0]^T = \nabla H \neq 0$ in the new coordinates.
 To study the Jacobian of $\Phi_{\Delta t}$ in the coordinates $(\theta,\phi)$ we let
\(
(\theta_{1}, \phi_{1}) = \Phi_{\Delta t}( \theta_0, \phi_0)
\)
for $(\theta_0,\phi_0) \in U$.
Then the Jacobian of $\Phi_{\Delta t}$ is
\[
\left[ \begin{array}{cc}
1 & 0 \\
\frac{\partial \phi_1}{\partial \theta_0} & \frac{\partial \phi_1}{\partial \phi_0}
\end{array}
\right].
\]
Area conservation implies that the determinant is
$1$ and so
\(
\partial \phi_1/\partial \phi_0 = 1
\)
 for all $\theta_0, \phi_0$.  This yields $\phi_1 = \phi_0 + \tau$ for all $\phi_0$ for some $\tau$ which may depend on $\theta_0$.
Hence, if we choose $c$ so that $c \Delta t/H'(\theta_0) = \tau$  we  have
\[
\Phi_{\Delta t}(\theta_0,\phi_0) = (\theta_0, \phi_0 + \tau) = S^{c \Delta t}(\theta_0,\phi_0),
\]
for all $(\theta_0,\phi_0) \in \Gamma$,
as required. $\Diamond$

For the main result we will use the following very weak definition of consistency.
\begin{defn} \label{defn:consistent}
Let $S^t$ be the flow map for the Hamiltonian $H: \mathbb{R}^{2n} \rightarrow \mathbb{R}$.
An integrator $\Phi$ is consistent at $(q,p) \in \mathbb{R}^{2n}$ if  
\[
\lim_{\Delta t \rightarrow 0} \frac{ \|S^{\Delta t}(q,p) - \Phi(H,\Delta t,(q,p)) \|}{\Delta t}  = 0.
\]
\end{defn}

\begin{thm} \label{thm:main}
Let $\Phi$ be an integrator for a Hamiltonian system in $\mathbb{R}^2$.
Let $H(q,p)=p^2/2$.
Let $U$ be an open neighbourhood of $\{(q,1)| q \in \mathbb{R} \}$.
Suppose:
\begin{enumerate}
\item $\Phi$ is consistent at $(0,1)$ (Definition \ref{defn:consistent}).
\item $\Phi$ is computationally reasonable (Definition \ref{defn:reasonable}).
\item $\Phi$ conserves energy for any $H$, $(q,p)$, and $\Delta t$ for which it is defined.
\item  
For  $\Delta t \leq \Delta t_0$,  $\Phi(H,(q,p),\Delta t)$ is defined, depends continuously on $(q,p)$, and conserves volume on $U$.  
\end{enumerate}
Then for all sufficiently small $\Delta t$ there is a 
$C_0^\infty$ function $V: \mathbb{R} \rightarrow \mathbb{R}$ such
that if $\widetilde{H}(q,p)=p^2/2 + V(q)$ then $\Phi(\widetilde{H},\cdot,\Delta t)$ is not simultaneously defined, continuous,
and volume-conserving on $U$.   For each such $\Delta t$, the $V$ constructed can be
replaced by $\lambda V$ for $\lambda \in (0,1)$ and the same result holds.
\end{thm}

{\bf Proof: }
For any $\Delta t \in (0, \Delta t_0]$, since $\Phi(H,\cdot,\Delta t)$ conserves energy and is continuous on $U$, the
set $\{(q,1) | q\in \mathbb{R} \}$ is mapped onto itself. 
As $\Phi(H,\cdot,\Delta t)$ conserves volume, Lemma~\ref{lem:reparam} shows that  it is identical to the flow of the original Hamiltonian system on $\{(q,1)| q\in \mathbb{R}\}$ with a rescaling of time:
\[
\Phi(H,(q,1),\Delta t) = (q+c \Delta t, 1),
\]
for all $q \in \mathbb{R}$,
where $c$ does not depend $q$.  
The consistency condition at $(0,1)$ implies that for small enough $\Delta t$ we have that $c>0$.  From now on, we assume $\Delta t$ is small enough so that $c>0$.

Consider the integrator applied to $H$ at the point $(q,p)=(0,1)$.  
Since $\Phi$ is computationally reasonable (Definition~\ref{defn:reasonable}), there are a
finite number of points $(r^j,s^j) \in \mathbb{R}^2$, $j=1,\ldots,m$,  such
that $\Phi$ only depends on $H$ at these points.
Choose a $q_0 \in \mathbb{R}$ big enough so that the interval  $[q_0,q_0+ c
  \Delta t ]$ contains none of the points $r^j, j=1,\ldots,m$ and is
disjoint from the interval $[0,c \Delta t]$. 

Consider the integrator applied to $H$ at the point $(q_0,1)$. There
are points $(\bar{r}^j,\bar{s}^j), j=1,\ldots,\bar{m}$  such
that $\Phi(H,(q_0,1),\Delta t)$ only depends on $H$ at these points.  
Let $V$ be a $C^\infty_0$ function such that
\begin{enumerate}
\item $V(q)=0$ for $q$ not in $[q_0, q_0+c \Delta t]$, and $V(\bar{r}^j)=0, j=1,\ldots,\bar{m}$,
\item $0 \leq V(q) < 1/2$ for all $q$,
\item for some $\epsilon>0$, $V(q)>0$ for $q \in (q_0,q_0+\epsilon)$.
\end{enumerate}
Note that multiplying $V$ by any factor in $(0,1)$ gives a function
satisfying the same conditions.  

Let $\widetilde{H}(q,p) = p^2/2 + V(q)$ and 
let $\widetilde{S}^t$ denote the flow map of the system with this Hamiltonian.
Now if $\Phi(\widetilde{H},\cdot,\Delta t)$ were defined, continuous
and volume-conserving on $U$, then by 
Lemma~\ref{lem:reparam} 
\[
\Phi(\widetilde{H},(q,1),\Delta t) = \widetilde{S}^{d \Delta t} (q,1)
\]
for all $q\in \mathbb{R}$ for some constant $d$.  We will show that
this is impossible. 

First note that the integrator gives the same result for both
Hamiltonians at $(0,1)$.  This is because here the result of the
integrator only depends on $H$ only at the points $(r^j,s^j)$, at which
$H$ and $\widetilde{H}$ agree. 
Since
\[
(d \Delta t,1) = \Phi(\widetilde{H},(0,1),\Delta t) = \Phi(H,(0,1),\Delta t) = (c \Delta t,1)
\]
we must have $c=d$.

On the other hand, the integrators also give the same results for both Hamiltonians at $(q_0,1)$.   This implies $c \neq d$ in the following way.  We have that 
\begin{eqnarray*} (q_0+c \Delta t,1) & = & \Phi (H,(q_0,1),\Delta t) \\
 & = & \Phi(\widetilde{H},(q_0,1),\Delta t)\\
 & = & \widetilde{S}^{d \Delta t} (q_0,1) \\ 
 & \neq & S^{d \Delta t}(q_0,1)\\ 
 & = & (q_0+d \Delta t,1).
\end{eqnarray*}
The inequality follows from Lemma~\ref{lem:differ}, since the
Hamiltonians differ on the interval $[q_0,q_0 + \epsilon]$. 
So $c \neq d$.  This is a contradiction.  Therefore,
$\Phi(\widetilde{H},\cdot,\Delta t)$ cannot be simultaneously defined, continuous, and
volume-conserving on $U$. $\Diamond$


The following lemma asserts the intuitively clear fact that if $V(q)$ is positive for $q \in [q_0, q_0+\epsilon]$ then the flows starting from $(q_0,1)$ corresponding to $H$ and $\widetilde{H}$ are different for positive times.

\begin{lem} \label{lem:differ}
Let $V: \mathbb{R}^{2n} \rightarrow \mathbb{R}$ be $C^\infty_0$ with $V(q_0)=0$, $1/2>V(q) \geq 0$ for all $q\in
\mathbb{R}$ and $V(q)>0$ for all $q \in [q_0, q_0+\epsilon]$. 
Let $H(q,p) = p^2/2$ and $\widetilde{H}(q,p) = p^2/2 + V(q)$.  Let
$S^t$ and $\widetilde{S}^t$ be the respective Hamiltonian flow maps of $H$ and $\widetilde{H}$.  Then
$S^t(q_0,1) \neq \widetilde{S}^t(q_0,1)$ for $t> 0$. 
\end{lem}

{\bf Proof.}
The trajectory for the Hamiltonian $H$ as a function of time is $(q(t),p(t))=(t,1)$ for all $t$.  
Letting $\tilde{q}(t)$ describe the position for the Hamiltonian $\widetilde{H}$, the usual solution technique gives 
\[
t = \int_{q_0}^{\tilde{q}} \frac{1}{\sqrt{1-2V(x)}} dx=: \widetilde{F}(\tilde{q})
\]
The function $\widetilde{F}$ is strictly increasing and so has a well defined inverse.
We can write
\(
\tilde{q}(t) = \widetilde{F}^{-1}(t).
\)
Now $\widetilde{F}'(\tilde{q}) > 1$ for $\tilde{q} >0$, so $\tilde{q}'(t) < 1$ for $t>0$.  Hence $\tilde{q}(t) < t$ for $t>0$ and the two flow maps cannot be equal for $t>0$.
$\Diamond$

\section{Multiple Degrees of Freedom}
\label{sec:multi}

In the previous section we showed that there can be no general energy-
and volume-conserving integration schemes because there are no integrators that conserve energy and volume for all Hamiltonians in $\mathbb{R}^2$.
However, suppose we ask if such integrators exists for Hamiltonian systems
of dimension $2n$, $n>1$.  We conjecture that
a result like Theorem~\ref{thm:main} still holds in this case.
However, the method of proof for the $n=1$ case does not
extend to this case. 

Instead we will state two conditions on an integrator, either one of
 which prevents it from being volume- and energy-conserving for
 general Hamiltonian systems in $\mathbb{R}^{2n}$.
Both of these conditions are desirable for an integrator to have, but
 unlike computational reasonibility, it is not difficult to imagine a practical
 integrator that did not satisfy them.
The proof of the theorems in this section will work by showing that,
 if a computational reasonable energy-conserving
 integrator with either condition exists, a
 special Hamiltonian in $\mathbb{R}^{2n}$ can be constructed for which
 it does not conserve volume.

In the following, fix $n \geq 2$.  Let $q,p \in \mathbb{R}^{n}$ have
components 
$q_i, p_i \in \mathbb{R}$.  Define the projections $\pi_i$, $i=1,\ldots,n$ by 
$\pi_i(q,p)=(q_i,p_i)$.

The first condition asserts that if variables $(q_i,p_i)$ do not occur in the Hamiltonian function then the integrator does not change their value.

\begin{cond} \label{cond:untouch}
If $H(q,p)$ is independent of the variables $(q_i,p_i)$
then 
\[
\pi_i \Phi( H, (q,p), \Delta t) = (q_i,p_i)
\]
for all $(q,p)$ and $\Delta t$ for which $\Phi$ is defined. 
\end{cond}

\begin{thm}
Let $\Phi$ be an integrator for Hamiltonian systems in $\mathbb{R}^{2n}$. Let $H(q,p)=p_1^2/2$.   Let $U$ be an open neighbourhood of the set 
\[
\left\{ \left( (q_1,0, \ldots,0), (1,0,\ldots,0) \right) | q_1 \in \mathbb{R} \right\}.
\]
Suppose
\begin{enumerate}
\item
$\Phi$ is consistent (Definition 3).
\item 
$\Phi$ is computationally reasonable (Definition 2).
\item
$\Phi$ conserves energy for any $H$, $(q,p)$ and $\Delta t$ for which it is defined.
\item
$\Phi$ satisfies Condition~\ref{cond:untouch}
\item
For sufficiently small $\Delta t$, $\Phi(H,(q,p),\Delta t)$ is defined, continuous, and conserves volume for $(q,p) \in U$.  
\end{enumerate}
Then for sufficiently small $\Delta t$, there is a $C_0^\infty$ function $V: \mathbb{R} \rightarrow \mathbb{R}$ such that if $\widetilde{H}(q,p) = p_1^2 + V(q_1)$, then $\Phi(\widetilde{H},\cdot, \Delta t)$ is not simultaneously defined, continuous, and volume-conserving on $U$.  For each such $\Delta t$, the $V$ constructed can be replaced by $\lambda V$ for $\lambda \in (0,1)$ and the same result holds. \end{thm}

{\bf Proof: } We will use the hypothesized integrator on $\mathbb{R}^{2n}$ to construct an integrator $\phi$ on $\mathbb{R}^2$ satisfying the conditions of Theorem~\ref{thm:main}.
Let $H:\mathbb{R}^2 \rightarrow \mathbb{R}$ be a given Hamiltonian function.  
Define $H^*: \mathbb{R}^{2n} \rightarrow \mathbb{R}$ by $H^*(q,p)= H(q_1,p_1)$.
We define the integrator $\phi$ by
\[
\phi(H,(q_1,p_1),\Delta t) = \pi_1 \Phi(H^*,((q_1, 0 ,\ldots,0),(p_1,0,\ldots,0)),\Delta t).
\]

Let $\bar{U} \subset \mathbb{R}^2$ be given by
\[
\bar{U} := \left\{ (q_1,p_1) | ((q_1,0,\ldots,0),(p_1,0,\ldots,0)) \in U \right\}.
\]
It is straightforward to check that $\phi$ satisfies the conditions of Theorem~\ref{thm:main} on $\bar{U}$.  Thus, by the theorem,
we have an arbitrarily small $C^\infty_0$ function $V$ such that $\phi$ is not simultaneously defined, continuous, and volume-conserving on $\bar{U}$ for $\widetilde{H}(q_1,p_1)=p_1^2+V(q_1)$.    

Let $\widetilde{H}^*(q,p) = p_1^2/2 + V(q_1)$ for $(q,p)\in \mathbb{R}^{2n}$.
Now suppose that $\Phi(\widetilde{H}^*,\cdot,\Delta t)$ is defined, continuous, and volume preserving on $U$. We will derive a contradiction by showing this implies that $\phi(\widetilde{H},\cdot,\Delta t)$ is, in fact, defined, continuous, and volume-conserving on $\bar{U}$.

First note that $\Phi(\widetilde{H}^*,\cdot,\Delta t)$ being defined and continuous on $U$ implies that $\phi(\widetilde{H},\cdot,\Delta t)$ is defined and continuous on $\bar{U}$.
To check volume conservation, note that the Jacobian of the map $\Phi(\widetilde{H}^*,\cdot,\Delta t)$ has structure
\[
\left[ 
\begin{array}{cc}
J_{11} & J_{12} \\
0        & I
\end{array}
\right]
\]
where we have put the variables in order $(q_1,p_1,\ldots,q_n,p_n)$ and $J_{11}$ is a 2-by-2 matrix.  
Since the determinant of this matrix is 1 by volume-conservation, the determinant of $J_{11}$ must be 1.  But $J_{11}$ is the Jacobian of $\phi(\widetilde{H},\cdot,\Delta t)$, so this latter map must be area preserving.  This contradicts our earlier assumption. $\diamond$

 The second condition states that if the
Hamiltonian system consists of $n$ identical uncoupled one-degree-of-freedom
systems, then the integrator itself should consist of $n$ identical uncoupled
maps on the state-space of each subsystem. 

\begin{cond}  \label{cond:prod}
If $H(q,p)= \sum_{i=1}^n h(q_i,p_i)$ 
then there is an integrator $\phi$ on $\mathbb{R}^2$ such that
\[
\pi_i  \Phi(H,(q,p),\Delta t) = \phi(h,(q_i,p_i),\Delta t)
\]
for $i=1,\ldots,n$, for any $(q,p)$ and $\Delta t$ for which the integrator is defined.
\end{cond}

Though this is certainly a nice property for the integrator to have
(since the flow map has the same property) there are many integrators for
which it does not hold.  For example, step-and-project methods may
not satisfy this condition, even if the underlying one-step method does.


\begin{thm}
Let $\Phi$ be an integrator for Hamiltonian systems in $\mathbb{R}^{2n}$. Let $H(q,p)=\sum_i p_i^2/2$.   Let $U \subset \mathbb{R}^{2n}$ be an open neighbourhood of the set
\[
\left\{ \left((q_1, \ldots,q_1), (1,\ldots,1)\right) | q_1 \in \mathbb{R} \right\}.
\]
Suppose
\begin{enumerate}
\item
$\Phi$ is consistent (Definition 3).
\item 
$\Phi$ is computationally reasonable (Definition 2).
\item
$\Phi$ conserves energy for any $H$, $(q,p)$ and $\Delta t$ for which it is defined.
\item
$\Phi$ satisfies Condition~\ref{cond:prod}
\item
For sufficiently small $\Delta t$, $\Phi(H,(q,p),\Delta t)$ is defined, continuous, and conserves volume for $(q,p) \in U$.  
\end{enumerate}
Then for sufficiently small $\Delta t$, there is a $C_0^\infty$ function $V: \mathbb{R} \rightarrow \mathbb{R}$ such that if $\widetilde{H}(q,p) = \sum_{i=1}^n (p_i^2 + V(q_i))$ then $\Phi(\widetilde{H},\cdot, \Delta t)$ is not simultaneously defined, continuous, and volume-conserving on $U$.  For each such $\Delta t$, the $V$ constructed can be replaced by $\lambda V$ for $\lambda \in (0,1)$ and the same result holds.
\end{thm}

{\bf Proof: } This theorem is proven analogously to the previous theorem.  For any $H: \mathbb{R}^{2} \rightarrow \mathbb{R}$ we define $H^* :\mathbb{R}^{2n} \rightarrow \mathbb{R}$ by $H^*(q,p)= \sum_i H(q_i,p_i)$.
We define the integrator $\phi$ by
\[
\phi(H,(q_1,p_1),\Delta t) = \pi_1 \Phi(H^*,((q_1,\ldots,q_1),(p_1,\ldots,p_1)),\Delta t).
\]
We define $\bar{U} \in \mathbb{R}^2$ by
\[
\bar{U} := \left\{ (q_1,p_1) | \left( (q_1,\ldots,q_1), (p_1,\ldots,p_1)\right) \in U \right\}.
\]
As in the proof of the previous theorem,
$\phi$ and $\bar{U}$ satisfy the conditions of Theorem~\ref{thm:main}.  Thus, by the theorem,
we have an arbitrarily small $C^\infty_0$ function $V$ such that $\phi$ is not defined continuous and volume-conserving on $\bar{U}$ for $\widetilde{H}=p_1^2+V(q_1)$.    

Let $\widetilde{H}^*(q,p) = \sum_i (p_i^2/2 + V(q_i))$ for $(q,p)\in \mathbb{R}^{2n}$.
Now suppose that $\Phi(\widetilde{H}^*,\cdot,\Delta t)$ is defined, continuous, and volume preserving on $U$. We will derive a contradiction.

Now $\phi(\widetilde{H},\cdot, \Delta t)$ is defined and continuous on $\bar{U}$.
To check volume conservation, note that the Jacobian of the map $\Phi(\widetilde{H}^*,\cdot,\Delta t)$ in this case has structure
\[
\left[ 
\begin{array}{cccc}
J_{11} & 0 &  \ldots & 0  \\
0         & J_{22} & \ldots & 0 \\
\vdots &  \vdots & \ddots & 0 \\
0         &  \ldots  & 0   & J_{nn}  
\end{array}
\right]
\]
where we have put the variables in order $(q_1,p_1,\ldots,q_n,p_n)$ and each $J_{ii}$ is 2-by-2.  Since the determinant of this matrix must be 1 by volume-conservation and the determinants of the $J_{ii}$ are identical, the determinant of $J_{11}$ must be $\pm 1$.  
As in the proof of the previous theorem this implies $\phi(\widetilde{H},\cdot,\Delta t)$ is area conserving on $\bar{U}$ which is a contradiction.
 $\diamond$

{\bf Acknowledgments.}  The author thanks Nilima Nigam for her comments on this work.  The author was supported by an NSERC Discovery Grant.

\end{document}